\documentclass[a4paper,12pt]{article}

\pdfoutput=1

\usepackage[paper=letterpaper,margin=1in]{geometry}

\usepackage{amsmath,amssymb,amsfonts,epsfig,cite,setspace,bigstrut,longtable,array,breqn,color}
\usepackage{todonotes}

\usepackage{hyperref,caption,setspace}

\usepackage{cancel}
\usepackage{pdfpages,lipsum}

\begin{document}



\vskip 0.25in

\newcommand{\nn}{\nonumber}
\newcommand{\tr}{\mathop{\rm Tr}}
\newcommand{\comment}[1]{}

\newcommand{\cM}{{\cal M}}
\newcommand{\cW}{{\cal W}}
\newcommand{\cN}{{\cal N}}
\newcommand{\cH}{{\cal H}}
\newcommand{\cK}{{\cal K}}
\newcommand{\cZ}{{\cal Z}}
\newcommand{\cO}{{\cal O}}
\newcommand{\cA}{{\cal A}}
\newcommand{\cB}{{\cal B}}
\newcommand{\cC}{{\cal C}}
\newcommand{\cD}{{\cal D}}
\newcommand{\cT}{{\cal T}}
\newcommand{\cV}{{\cal V}}
\newcommand{\cE}{{\cal E}}
\newcommand{\cF}{{\cal F}}
\newcommand{\cX}{{\cal X}}
\newcommand{\IA}{\mathbb{A}}
\newcommand{\IP}{\mathbb{P}}
\newcommand{\IQ}{\mathbb{Q}}
\newcommand{\IH}{\mathbb{H}}
\newcommand{\IR}{\mathbb{R}}
\newcommand{\IC}{\mathbb{C}}
\newcommand{\IF}{\mathbb{F}}
\newcommand{\IV}{\mathbb{V}}
\newcommand{\II}{\mathbb{I}}
\newcommand{\IZ}{\mathbb{Z}}
\newcommand{\re}{{\rm~Re}}
\newcommand{\im}{{\rm~Im}}

\let\oldthebibliography=\thebibliography
\let\endoldthebibliography=\endthebibliography
\renewenvironment{thebibliography}[1]{%
\begin{oldthebibliography}{#1}%
\setlength{\parskip}{0ex}%
\setlength{\itemsep}{0ex}%
}%
{%
\end{oldthebibliography}%
}

\newtheorem{theorem}{\bf THEOREM}
\def\thetheorem{\thesection.\arabic{theorem}}
\newtheorem{proposition}{\bf PROPOSITION}
\def\thetheorem{\thesection.\arabic{proposition}}
\newtheorem{observation}{\bf OBSERVATION}
\def\thetheorem{\thesection.\arabic{observation}}
\newtheorem{definition}{\bf DEFINITION} 
\def\thetheorem{\thesection.\arabic{DEFINITION}}
\newtheorem{example}{\bf EXAMPLE} 
\def\thetheorem{\thesection.\arabic{EXAMPLE}}
\newtheorem{conjecture}{\bf CONJECTURE} 
\def\thetheorem{\thesection.\arabic{CONJECTURE}}
\newtheorem{remark}{\bf REMARK} 
\def\thetheorem{\thesection.\arabic{REMARK}}

\def\theequation{\thesection.\arabic{equation}}
\newcommand{\setall}{\setcounter{equation}{0}
        \setcounter{theorem}{0}}
\newcommand{\setequation}{\setcounter{equation}{0}}
\renewcommand{\thefootnote}{\fnsymbol{footnote}}

\newcommand{\seteq}{\mathbin{:=}}
\newcommand{\GL}{\operatorname{GL}}
\newcommand{\Sp}{\operatorname{Sp}}
\newcommand{\USp}{\operatorname{USp}}
\newcommand{\GSp}{\operatorname{GSp}}
\newcommand{\U}{\operatorname{U}}
\newcommand{\SU}{\operatorname{SU}}
\newcommand{\SO}{\operatorname{SO}}
\newcommand{\End}{\operatorname{End}}

\begin{titlepage}

~\\
\vskip 1cm

\begin{center}
{\Large \bf Machine-Learning the Sato--Tate Conjecture}
\end{center}
\medskip

\renewcommand{\arraystretch}{0.5}

\vspace{.4cm}
\centerline{
{\large Yang-Hui He, Kyu-Hwan Lee, Thomas Oliver}
}
\vspace*{3.0ex}

\vspace{10mm}

\begin{abstract}
We apply some of the latest techniques from machine-learning to the arithmetic of hyperelliptic curves. More precisely we show that, with impressive accuracy and confidence (between $99$ and $100$ percent precision), and in very short time (matter of seconds on an ordinary laptop), a Bayesian classifier can distinguish between Sato--Tate groups given a small number of Euler factors for the $L$-function. Our observations are in keeping with the Sato-Tate conjecture for curves of low genus.  For elliptic curves, this amounts to distinguishing generic curves (with Sato--Tate group $\SU(2)$) from those with complex multiplication. In genus $2$, a principal component analysis is observed to separate the generic Sato--Tate group $\USp(4)$ from the non-generic groups. Furthermore in this case, for which there are many more non-generic possibilities than in the case of elliptic curves, we demonstrate an accurate characterisation of several Sato--Tate groups with the same identity component. Throughout, our observations are verified using known results from the literature and the data available in the LMFDB. The results in this paper suggest that a machine can be trained to learn the Sato--Tate distributions and may be able to classify curves much more efficiently than the methods available in the literature. 
\end{abstract}

\end{titlepage}

\tableofcontents

\section{Introduction \& Summary}

There is a strong tradition of machine aided computation in number theory, which has been used to formulate and verify a wide range of arithmetic conjectures. In this paper, we pursue a data-driven approach to a classification problem in arithmetic geometry. In particular, we demonstrate that a Bayesian classifier can efficiently and accurately distinguish Sato--Tate groups of genus $1$ and $2$ curves.

The original Sato--Tate conjecture is concerned with the distribution of Euler factors associated to elliptic curves over number fields. In recent years, there has been remarkable progress made towards this conjecture, which would be a corollary to establishing certain analytic properties of symmetric power $L$-functions. The necessary analytic behaviour would be a consequence of Langlands functoriality. In fact, it is sufficient to prove potential automorphy (automorphy after base change to a field extension). This idea has been used to establish the Sato--Tate conjecture for elliptic curves over various fields \cite{Tay08}, \cite{HSBT},  \cite{ACCGH+}. There is also a body of literature for more general Hilbert modular forms.

A precise analogue of the Sato--Tate conjecture for genus 2 curves was formulated in \cite{KS09,FKRS}. In this context, there are 52 possible distributions corresponding to various endomorphism types of the Jacobian. For genus 2 curves defined over $\mathbb{Q}$, the number of possibilities is reduced to 34. Each distribution can be described by the Haar measure of a compact Lie group known as the Sato--Tate group. The generalized Sato--Tate conjecture asserts that the distribution of the Euler factors converges to the distribution of the characteristic polynomials of random matrices in the Sato--Tate group. 

As with elliptic curves, the Sato--Tate conjecture for genus 2 curves would follow from the Langlands functoriality conjectures \cite[Section~1.7]{FKRS}. The Sato--Tate conjecture for non-generic genus 2 curves over $\mathbb Q$ has been established by C. Johansson and N. Taylor \cite{Joh17, Tay20}. Conditional on the Sato--Tate conjecture, one may compute the Sato--Tate group of a genus 2 curve by evaluating moments of the coefficients appearing in normalized Euler factors and comparing to the corresponding statistics for characteristic polynomials of random matrices. This approach was adopted in \cite{KS09}. The Sato--Tate groups on the LMFDB were confirmed by an unconditional approach in \cite{CMSV} to compute the real endomorphism algebra. See \cite[Section~4.4]{BSSVY} for more explanation.

In parallel to the above developments, a recent programme of machine-learning mathematical structures was initiated in \cite{He:2017aed,He:2017set}. Whilst this was originally motivated by computing topological invariants of Calabi--Yau compactifications in superstring theory \cite{He:2017set,Krefl:2017yox,Ruehle:2017mzq,Carifio:2017bov} (q.v., \cite{He:2018jtw} for a summary), the idea of using machine-learning for pattern-recognition and conjecture-raising has been applied to various branches of mathematics, such as representation theory \cite{He:2019nzx}, graph theory \cite{He:2020fdg}, metric geometry \cite{Ashmore:2019wzb}, knot invariants \cite{Jejjala:2019kio}, quiver mutations \cite{Bao:2020nbi}, etc.
The reader is also pointed to interesting early \cite{shanker} and recent \cite{KV} experiments in neural-network explorations of the famous zeros of the Riemann zeta function. Machine learning techniques were applied to databases of elliptic curves in \cite{Alessandretti:2019jbs}. In that work, the data consisted of the Weierstra\ss\ coefficients for each curve. These coefficients vary in size dramatically, which partially accounted for the difficulty in mining the data.

In this paper, we study the (conditional) computation of Sato--Tate groups via machine learning techniques. Naturally, this approach requires a large amount of data to train the algorithm. Much data can be sourced from the LMFDB, which enables a classifier to efficiently distinguish curves belonging to certain pairs of Sato--Tate groups \cite{lmfdb}. 
There are not enough examples of curves for the other Sato--Tate groups for a full classification, and so we turn to random matrices to generate our training data. Using this, we are able to establish a finer classification which, for example, can distinguish curves from 5 Sato--Tate groups with the same identity component. Applying the same method, we can train a classifier with data coming from random matrices of the 34 Sato--Tate groups for genus 2 curves over $\mathbb Q$. Nevertheless, for the present, we are unable to verify the accuracy of a full 34-way classification due to a lack of available data.

The organization of this paper is as follows. In Section \ref{sec-background}, the generalized Sato--Tate conjecture for arithmetic curves will be reviewed as the main mathematical background for this paper. In Section \ref{sec-LMFDB}, machine-learning techniques will be applied to certain binary classifications of curves. In Section \ref{sec-finer}, we go beyond the binary classification and consider a multi-way classification of genus 2 curves corresponding to Sato--Tate groups with a common identity component. Throughout, we compare the machine learning method with other approaches to computation of the Sato--Tate groups. The general theme is that machine-learning requires a significantly smaller input to determine the Sato--Tate group. This is of course after the classifier has been trained, which takes only a matter of seconds on an ordinary laptop.

\subsection*{Acknowledgements}
We thank \'Alvaro Lozano-Robledo, Andrew Sutherland and Chris Wuthrich for helpful discussions and useful comments.
YHH is indebted to STFC UK, for grant ST/J00037X/1, KHL is partially supported by a grant from the Simons Foundation (\#712100), and TO acknowledges support from the EPSRC through research grant EP/S032460/1.

\section{Background} \label{sec-background}

In this section we review the essential mathematical theory which constitutes the main theme of this paper.

\subsection{CM elliptic curves}\label{sec:cmcurves}

Let $\mathcal{E}$ be an elliptic curve over $\mathbb{Q}$. With minor modifications, it is possible to replace $\mathbb{Q}$ with any number field. Recall 
that the (Hasse--Weil) $L$-function depends only on the isogeny class of $\cE$ and captures many of its deep arithmetic properties. This function is given by an Euler product:
\begin{equation}\label{eq.EllipticEuler}
L(s, \cE) 
=\prod\limits_{p \mid N} (1 - a_p p^{-s})^{-1}\prod\limits_{p \nmid N} (1 - a_p p^{-s} + p^{-2s})^{-1} \ ,
\end{equation}
where $N$ is the conductor, which controls primes of good and bad reduction.

The elliptic curve $\cE$ is said to have CM if its ring of endomorphisms is strictly larger than the ring of integers. In terms of the Sato--Tate conjecture, a CM elliptic curve has the distribution of normalized Euler factors converging to that of characteristic polynomials of random matrices in the normalizer $N(\U(1))$ of $\U(1)$ in $\SU(2)$, while a non-CM curve has the distribution of normalized Euler factors converging to that of characteristic polynomials of $\SU(2)$. In a rigorous sense explained in the next subsection, non-CM curves are generic, while CM curves are exceptional.
 
An elliptic curve $\cE$ has CM, or equivalently, its Sato--Tate group is $N(\U(1))$ if its $j$-invariant is one of 13 integers listed 
\footnote{
Namely:
\[
\begin{array}{ll}
& 0,\quad 2^4 \ 3^3\  5^3, \quad  -2^{15 }\ 3 \ 5^3, \quad 2^6 \ 3^3, \quad 
 2^3 \ 3^3 \ 11^3, \quad  -3^3 \ 5^3, \quad 3^3 \ 5^3 \ 17^3, \quad  
 2^6 \ 5^3,\\ 
 & -2^{15}, \quad
  -2^{15} \ 3^3,\quad -2^{18} \ 3^3 \ 5^3,\quad 
-2^{15} \ 3^3 \ 5^3 \ 11^3, \quad  -2^{18} \ 3^3 \ 5^3 \ 23^3 \ 29^3  \ . \\
\end{array}
\]
}
in \cite[Appendix~A,~Section~3]{silverman2}.
This criterion for CM curves is based on the result of Heegner--Baker--Stark. The $j$-invariant is an elementary function in terms of the Weierstra\ss\ coefficients. On the other hand, the Dirichlet coefficients $\{ a_p \}$ encode CM in other ways. If $\cE$ over $\mathbb{Q}$ has CM by the integers in an imaginary quadratic number field $K$, then there is a Hecke character $\psi$ on $\mathbb{A}_K^{\times}$ such that $L(s,\mathcal{E})=L(s,\psi)$ \cite[Theorem~10.5(b)]{silverman2}. It follows from the Chebotarev density theorem that the following set has density $1/2$ in the set of primes:
\[
\pi(\mathcal E)=\{p\text{ prime} \ : \ a_p=0\}.
\] 
On the other hand, if $\mathcal E$ does not have CM then 
$\pi(\mathcal E)$ has density $0$ in the primes \cite{Serre} though $\pi(\mathcal E)$ is still infinite as demonstrated by Elkies \cite{Elkies}. We note that it is in fact possible to distinguish CM from non-CM, or, equivalently, to determine whether the Sato--Tate group is $N(\U(1))$ or $\SU(2)$, given only finitely many $a_p$. There is a large body of literature concerned with questions of this nature, building on \cite{LO75}. Our approach in this paper is to use machine-learning techniques.

\subsection{Generalized Sato--Tate conjecture}

In this section, we briefly overview the generalized Sato--Tate conjecture, in particular, for genus 2 curves over $\mathbb Q$. More details can be found in \cite{KS09,FKRS}.

Let $\mathcal C$ be a smooth, projective, geometrically irreducible algebraic curve of genus $g$ defined over $\mathbb Q$. (The elliptic curves over $\mathbb Q$ lay in the subclass of $g=1$.) 
For each prime $p$ where $\mathcal C$ has good reduction,
we define the zeta
function by
\begin{equation}
Z(\mathcal C/\mathbb F_p; T) = \exp \left ( \sum_{k=1}^\infty N_kT^k/k \right ) , 
\end{equation}
where $N_k$ is the number of the points on $\mathcal C$ over $\mathbb F_{p^k}$. It is well-known that the zeta function can be written in the form
\begin{equation}\label{eq.genus2ZL}
Z(\mathcal C/\mathbb F_p; T) = \frac{L_p(T)}{(1-T)(1-pT)} , 
\end{equation}
where $L_p \in \mathbb Z[T]$ is a  polynomial of degree $2g$ with constant term $1$. 
 In particular, when $g=1$, we have $L_p(T)=1-a_pT+pT^2$ where $a_p$ appears in the Euler factor of the $L$-function in \eqref{eq.EllipticEuler}.
If we set $\bar L_p(T):=L_p(p^{-1/2}T)$, then  we obtain
\begin{equation}\label{akp}
 \bar{L}_p(T)=T^{2g}+a_{1,p}T^{2g-1}+ a_{2,p}T^{2g-2}+ \cdots + a_{2,p} T^2+a_{1,p}T+1 \ .
\end{equation}
We see that this normalization renders the $L$-function {\it palindromic}.

Let $P_{\mathcal C}(N)$ be the set of primes $p \le N$ for which the curve $\mathcal C$ has good reduction.
For $1\le k \le g$ and $m\ge 0$, define
\begin{equation} \label{eq:akmg}
a_k(m;g):=\lim\limits_{N \rightarrow \infty} \frac{1}{|P_{\mathcal C}(N)|} \sum\limits_{p \in P_{\mathcal C}(N)} (a_{k,p})^m.
\end{equation}
Thus the values $a_k(m;g)$, $m \ge 0$, are  the $m^{\mathrm{th}}$ moments of the distribution of $a_{k,p}$.

The generalized Sato--Tate conjecture predicts  that curves of fixed genus $g$ are classified into certain families and that  $a_k(m;g)$ are all the same for curves in each family. In particular, there is a generic family of curves for each genus $g$, which is characterized by the property that the Jacobians of its members have the trivial endomorphism ring $\mathbb Z$.
When $g=1$, the generic family exactly consists of non-CM elliptic curves.

The generalized Sato--Tate conjecture predicts that the distributions of $\bar L_p(T)$ are actually the same as the distributions of the characteristic polynomials of random matrices. To be precise, let us consider  the group $\USp(2g)$ with the Haar probability measure.  Let 
\begin{equation}\label{charpoly}
\det (I - x\gamma)=x^{2g}+c_{1}x^{2g-1}+ c_{2}x^{2g-2}+ \cdots + c_{2} x^2+c_{1}x+1
\end{equation}
be the characteristic polynomial of a random matrix $\gamma$ of $\USp(2g)$. For each $k=1,2, \dots, g$, let $X_k$ be the random variable corresponding to the coefficient $c_k$ and define $c_k(m;g)$ to be the $m^{\rm{th}}$ moment $\mathbf E[X_k^m]$, $m \in \mathbb Z_{\ge 0}$, of the random variable $X_k$. 

The following is the generalized Sato--Tate conjecture for the generic families.

\begin{conjecture}[\cite{KS}] \label{conj-ST}
Let $\mathcal C$ be a smooth projective curve of genus $g$.
Assume that $\mathcal C$ is in the generic family.
 Then, for each $k=1, 2, \dots, g$ and $m \ge 0$, we have
 \[ a_k(m;g)= c_k(m;g). \]
\end{conjecture}

In the case that $g=2$, a precise formula for $c_k(m;2)$ is given in \cite[Tables~9~\&~10]{FKRS}. Given a genus 2 curve, one may compute Euler factors for primes less than, say, $N$. The finite sum $\frac{1}{|P_{\mathcal C}(N)|} \sum\limits_{p \in P_{\mathcal C}(N)} (a_{k,p})^m$ provides an approximation to $a_k(m;2)$. Conditional on the Sato--Tate conjecture, we can check whether a curve is generic by comparison with the formula for $c_k(m;2)$. This identification is accurate up to a certain probability (discussed in Section~\ref{sec:Outlook}). We refer to this as the ``heuristic'' computation of the Sato--Tate group. 

\begin{example} \label{ex-ST}
The following genus $2$ curve (LMFDB label: 11109.a.766521.1) is from the generic family: 
\[
\mathcal C: y^2+(x^2+x)y=x^5-x^4+x^3-3x^2+2x-1.
\]
Conditional on the Sato--Tate conjecture the sequences $\{ a_k(m;2)\}$ are as follows:
\begin{align*}
a_1(m;2) :&\quad  1,0,1,0,3,0,14,0,84,0,594,0,4719, \dots \\
a_2(m;2):&\quad  1,1,2,4,10,27,82,268,940, \dots
\end{align*}
\end{example}

Aside from the generic family of curves whose distribution is (expected to be) given by $\USp(2g)$, there are exceptional families of curves. As mentioned in the previous subsection, the CM curves form the exceptional family when $g=1$, and the distribution is given by the normalizer $N(\U(1))$ of $\U(1)$ in $\SU(2) \cong  \USp(2)$.

For genus $2$ curves, there are a lot more of exceptional families. Kedlaya and Sutherland \cite{KS09} and later with Fit{\' e} and  Rotger \cite{FKRS} made a conjectural, exhaustive list of 34 compact subgroups of $\USp(4)$ that would classify all the distributions of Euler factors for genus 2 curves over $\mathbb Q$, and called the groups {\em Sato--Tate groups}. They determined the moment sequences $c_k(m;2)$, $k=1,2$, for each Sato--Tate group. In the process they investigated a huge number of genus 2 curves to heuristically observe that Euler factors have the same distributions as the Sato--Tate distributions, supporting their refined, generalized Sato--Tate conjecture. As with generic curves, one may heuristically compute the Sato--Tate group of any genus 2 curve by first computing an approximation to the moments and then comparing to the tables given in \cite[Tables~9~\&~10]{FKRS}. 

Since \cite{FKRS} appeared, the Sato--Tate conjecture for genus 2 curves over $\mathbb Q$ has been established by C. Johansson and N. Taylor \cite{Joh17, Tay20} except for the generic case $\USp(4)$. In particular, this means that the heuristic computation in these cases is no longer conditional (though it is still only valid up to a certain probability). The auto-correlation functions of the Sato--Tate distributions are computed in \cite{LO} using irreducible characters of symplectic groups, which provides an alternative way of characterizing the Sato--Tate distributions.

\begin{example}\label{eg.density}
In Section~\ref{sec:cmcurves}, we saw that non-generic elliptic curves were characterized by the density of vanishing coefficients. This can be predicted by computation of characteristic polynomials of cosets of the identity components as in \cite{LO}. For example, when $g=1$, the Sato--Tate group $N(\U(1))$ for CM curves has the coset decomposition 
\[ 
N(\U(1)) = \U(1) \sqcup J_2\U(1), 
\] 
where $J_2:=\scriptsize{\begin{pmatrix} 0&1 \\ -1&0 \end{pmatrix}}$, and the characteristic polynomial of the matrices from the coset $J_2 \U(1)$ is always $1+x^2$. This shows that $a_p=0$ with density $1/2$ for CM-curves. A similar analysis can be done for genus 2 curves by considering coset decompositions. 
\end{example}

In what follows, we define the Sato--Tate groups for genus 2 curves over $\mathbb Q$. 
We will adopt the same notations as in \cite{FKRS}. We take the group $\USp(4)$ to fix the symplectic form $\begin{pmatrix} 0 & I_2 \\ -I_2 & 0 \end{pmatrix}$, where $I_2$ is the $2 \times 2$ identity matrix. Let $E_{ij}$  be the $4 \times 4$ elementary matrix which has  $(i,j)$-entry equal to $1$ and other entries equal to $0$.
Set
\begin{align*}
\hat{h}_1&=E_{11}-E_{33}, &
\hat{h}_2&=E_{22}-E_{44} .
\end{align*}
We embed $\U(1)$ into $\USp(4)$ by
\begin{equation*} 
 u \longmapsto \mathrm{diag}(u,u,u^{-1},u^{-1}) .
\end{equation*}
For example, $e^{\pi i/n}$ is identified with  \[\mathrm{diag}(e^{\pi i/n},e^{\pi i/n},e^{-\pi i/n},e^{-\pi i/n}).\]
Embed  $\SU(2)$ and $\U(2)$ into $\USp(4)$ by
\begin{equation}\label{embed-su2} 
A \longmapsto \begin{pmatrix} A&0  \\ 0& \overline{A} \end{pmatrix},
\end{equation}
where $\overline A$ consists of the complex conjugates of the entries of $A$.

We fix an embedding \begin{equation} \label{embed-su2-su2} \SU(2) \times \SU(2) \hookrightarrow \USp(4)\end{equation} in such a way that the induced Lie algebra embedding $\mathfrak {sl}_2(\mathbb C) \times \mathfrak{sl}_2(\mathbb C) \rightarrow \mathfrak {sp}_4(\mathbb C)$ gives
\[ 
(h,0) \longmapsto \hat{h}_1 \quad \text{ and }\quad (0,h) \longmapsto \hat{h}_2 ,
\] 
where $h= \scriptsize{\begin{pmatrix} 1&0\\0&-1 \end{pmatrix}} \in \mathfrak{sl}_2(\mathbb C)$. From this, we also obtain the embeddings 
\[
\U(1)\times \SU(2) \hookrightarrow \USp(4) , \qquad \U(1)\times \U(1) \hookrightarrow \USp(4).
\]
Identify $\SU(2)$ with the group of unit quaternions via the isomorphism
\[  a+ b \, \pmb{\rm i} + c \, \pmb{\rm j} + d \, \pmb{\rm k}  \mapsto \begin{pmatrix} a+bi & c+di \\ -c+di & a-bi \end{pmatrix},  \quad a,b,c,d \in \mathbb R,   \]
and also identify them with the corresponding elements in $\USp(4)$ through the embedding $\SU(2) \hookrightarrow \USp(4)$ in \eqref{embed-su2}. For example, with this identification, we have
$ \mathbf j = \scriptsize { \begin{pmatrix} 0 & 1 & 0 & 0 \\ -1 & 0 & 0&0 \\ 0 &0&0&1\\ 0& 0&-1&0 \end{pmatrix} }$.
Set 
$
Q_1=\{ \pm 1, \pm \mathbf i, \pm \mathbf j, \pm \mathbf k, \tfrac 1 2 ( \pm 1 \pm \mathbf i \pm \mathbf j \pm \mathbf k \} )$ and
\begin{align*}
Q_2&=\left \{ \tfrac 1 {\sqrt 2} ( \pm 1 \pm \mathbf i), \tfrac 1 {\sqrt 2} ( \pm 1 \pm \mathbf j), \tfrac 1 {\sqrt 2} ( \pm 1 \pm \mathbf k),\tfrac 1 {\sqrt 2} ( \pm \mathbf i \pm \mathbf j),\tfrac 1 {\sqrt 2} ( \pm \mathbf i \pm \mathbf k),\tfrac 1 {\sqrt 2} ( \pm \mathbf j \pm \mathbf k) \right \}.
\end{align*}
 We write
$\pmb{\zeta}_{2n}= \begin{pmatrix} e^{\pi i/n} & 0 \\ 0 & e^{-\pi i/n} \end{pmatrix} \in \SU(2)$, and its embedded image in $\USp(4)$ will also be written as $\pmb{\zeta}_{2n}$.
Let
\begin{align*}  J& = \scriptsize { \begin{pmatrix} 0 & 0 & 0 & 1 \\ 0 & 0 & -1&0 \\ 0 &-1&0&0\\ 1& 0&0&0 \end{pmatrix} }, & \mathtt a &= {\scriptsize \begin{pmatrix}  0&0&1&0 \\  0&1&0&0 \\ -1&0&0&0 \\ 0&0&0&1 \end{pmatrix}},  & \mathtt b& = {\scriptsize \begin{pmatrix}  1&0&0&0 \\  0&0&0&1 \\ 0&0&1&0 \\ 0&-1&0&0 \end{pmatrix}},
& \mathtt c& = {\scriptsize \begin{pmatrix}  0&1&0&0 \\  -1&0&0&0 \\ 0&0&0&1 \\ 0&0&-1&0 \end{pmatrix}}.  \end{align*}

\begin{definition}[Sato--Tate groups] \label{def-ST}
With the notations above, the following table gives the definitions of the 34 Sato--Tate groups of genus 2 curves over $\mathbb Q$:
\begin{center}
\begin{tabular}{|rl|rl|}
\hline
$J(C_n)$ & $:= \langle \U(1), \pmb{\zeta}_{2n}, J \rangle $,  $n=2,4,6$ & $J(D_n)$ & $: = \langle J(C_n),  \mathbf j \rangle$, $n=2,3,4,6$\\ $J(T)$ & $:= \langle \U(1), Q_1, J \rangle$ & $J(O)$ & $:= \langle J(T), Q_2 \rangle$ \\ $ C_{n,1}$ & $ := \langle \U(1) , J \pmb{\zeta}_{2n} \rangle$, $n=2,6$ & $ D_{n,1}$ & $ := \langle \U(1), J \pmb{\zeta}_{2n},\,  \mathbf j \rangle$, $n =2,4,6$ \\
$ D_{n,2}$ & $ := \langle \U(1), \pmb{\zeta}_{2n}, J \mathbf j \rangle$, $n=3,4,6$ & $O_1$ & $ := \langle T, J Q_2 \rangle$ \\ $E_n$ & $:= \langle \SU(2), e^{\pi i/n} \rangle$, $n=1,2,3,4,6$ &  $J(E_n)$ & $ := \langle E_n, J \rangle $, $n=1,2,3,4,6$\\
$F_{\mathtt{a},\mathtt{b}}$ & $:= \langle \U(1) \times \U(1), \mathtt{a},\mathtt{b} \rangle $& $F_{\mathtt{ac}}$ & $:= \langle  \U(1) \times \U(1), \mathtt{ac} \rangle$ \\
$N(G_{1,3})$ & $:= \langle \U(1) \times \SU(2), \mathtt{a} \rangle $ & $ G_{3,3}$ & $:= \SU(2) \times \SU(2)$ \\ $N(G_{3,3})$ & $:= \langle G_{3,3}, J \rangle $ &  $\USp(4)$ & \\
\hline
\end{tabular}
\end{center} 
\end{definition}

\paragraph{Remark:}
We emphasize again that all the groups in the above table are subgroups of $\USp(4)$.
We will refer to the full $\USp(4)$ as the {\it generic} Sato-Tate group and the proper subgroups as the {\it non-generic}.

\newpage

\section{Distinguishing generic curves using the LMFDB} \label{sec-LMFDB}
In this section we describe a rudimentary binary classification using machine-learning techniques.

\subsection{Generic elliptic curves}

The latest LMFDB database has 3,064,705 elliptic curves over the rationals, 
which organize into 2,164,260 isogeny classes \cite[Elliptic curves over $\mathbb{Q}$]{lmfdb}. These curves are labeled by data of the form:
\begin{equation}
\{N, \ i, x\}
\end{equation}
where $N$ is the conductor, $i$ is a letter or double-letter designating the isogeny class, and $x$ is a number indexing the particular elliptic curve within the class (a typical entry, for instance, is `11a.1').  
For an elliptic curve, both its $L$-function (up to Euler factors at bad primes) and whether it has complex multiplication depend only on the isogeny class. Thus, for our present purpose, we will neglect the last numerical label $x$ and sometimes refer to the ``isogeny class of a curve'' simply as ``curve''. Of the some 2 million isogeny classes of elliptic curves in the database, only 2670 have CM: thus one can see that indeed this property is rather rare. 

Let us establish a dataset as follows. Take all primes up to 10,000 (there are 1229) and compute, using \cite{sage}, all coefficients $a_p$. Here we also include bad primes as their statistical impacts seem limited.
We then normalize the coefficients by $\widetilde{a_p} := a_p / \sqrt{p}$.
Now, take all 2670 curves with CM, and select with probability 0.001 from those without CM (which is therefore around  2300).
This gives a labeled dataset $\cD$ of around 5000 points:
\begin{equation}
\cD := \left\{
\left( \widetilde{a_p} \right)_{p < 10000} \longrightarrow \mbox{yes/no}
\right\}
\end{equation}
where yes/no refers to the simple binary category of having or not having CM.

Now, we can follow the standard steps of machine-learning (ML), which is to split $\cD$ into the disjoint union of a training set $\cT$ (taken a random sample) and a validation set $\cV$ (as the complement), and we take a 20-80 percent split:
\begin{equation}
\cD = \cT \sqcup \cV \ ,\qquad |\cT| = 20 \% |\cD| \ .
\end{equation}
The size of $\cV$ is large enough to check the validity of our results thoroughly. 
We tried a few architectures such as support vector machines and simple neural network classifiers, but found the best performance was achieved by a {\it Naive Bayes} classifier \footnote{
Interestingly, this is the same in the situation of machine recognition of cluster mutation \cite{Bao:2020nbi}.}.
We have also tried other standard classifiers, such as decision trees and nearest neighbours, but here, the Naive Bayes classifier performed best, and was able to achieve complete classification as we shall see shortly.
The reader is referred to \cite[Section~6.6.3]{Bayes} for detailed discussions and implementations of the algorithm.

We find that having seen 20\% of the $\widetilde{a_p}$-coefficients as lists of vectors, each of length 1229, and labeled accordingly as yes/no, the classifier, when validated on the remaining 80\%, achieves 100\% accuracy.
This is really the optimal situation. Ordinarily, a good classifier performs with precision (\% agreement) and confidence \footnote{
Matthew's phi-coefficient \cite{phi}, which essentially the square root of the chi-squared; the closer it is to 1, the better the fit, the closer it is to 0, the more random and ineffective the classification is.
We need to check this in addition to the naive precision in order to avoid false positives and false negatives.
}
in the 90's. 
But here, we consistently obtain 100\% accuracy with different random sampling of $\cT$. This suggests the ML algorithm has truly learned an underlying formula.
Moreover, the algorithm is performed using \cite{wolfram} on an ordinary laptop, in a matter of seconds.

\begin{figure}[!h!t!b]
\centerline{
\includegraphics[trim=12mm 0mm 0mm 0mm, clip, width=4in]{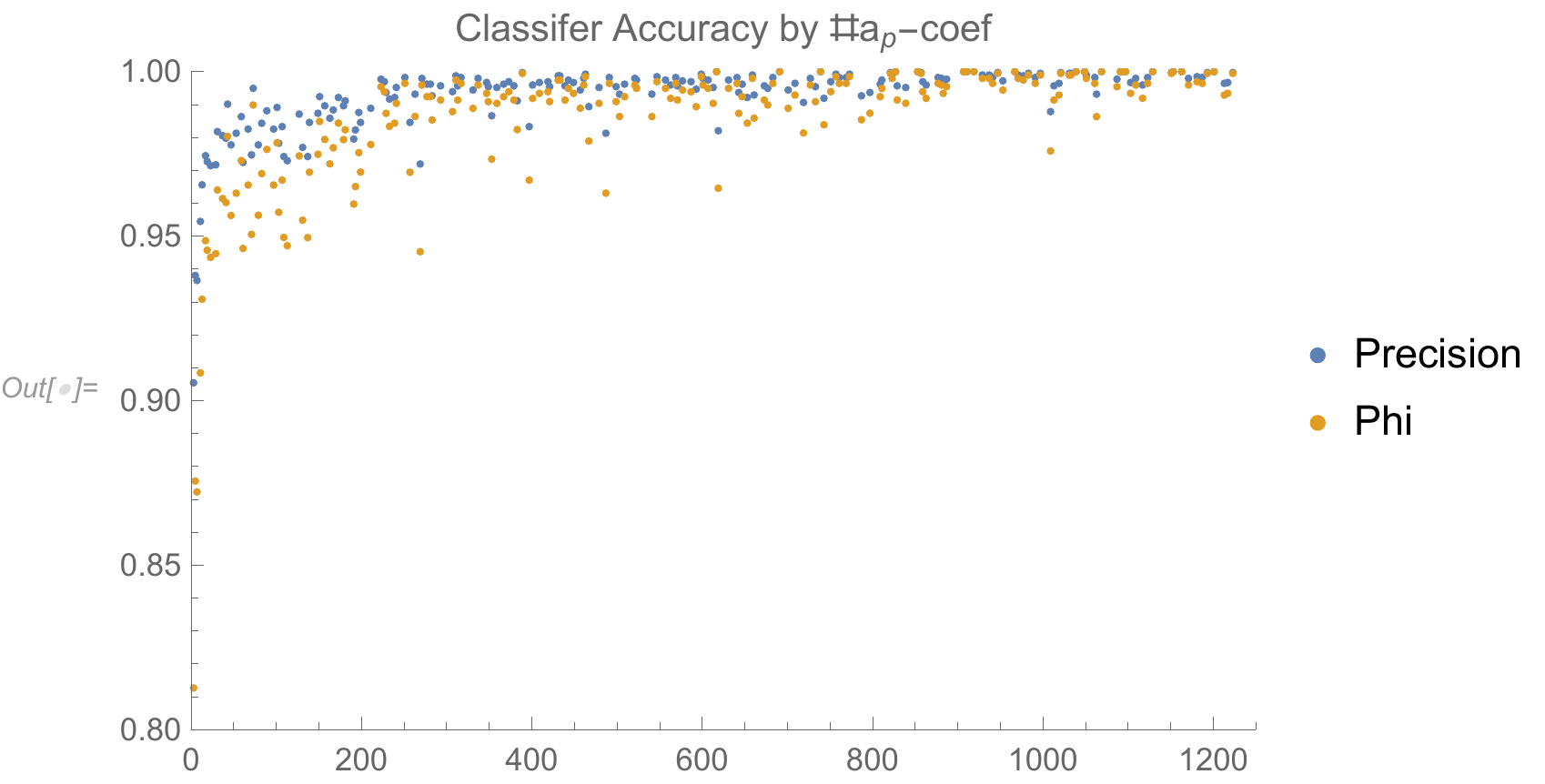}
}
\caption{{\sf {\small
The precision and confidence of the Naive Bayes classifier for the precision and confidence (Matthew's phi-coefficient) against the number of $a_p$ coefficients, for $p$ up to the value on the X-axis, for the elliptic curve seen.
}}
\label{f:elliptic-ap}}
\end{figure}

To get an idea of the learning, let us ask how the accuracies improve with increasing number of coefficients $a_p$ being presented to the training.
This is shown in Figure \ref{f:elliptic-ap}.
In other words, let us repeat the above Bayes classifier for truncated input data: instead of using all primes up to 10,000, we use up to the first 200 primes, in increments. While in the beginning the precision and phi are both low and sporadic, by the time we are training on primes up to 200 (i.e., only around 40 $a_p$ coefficients), we have stabilized to $>0.99$ accuracies. 

\subsection{Generic genus 2 curves}
Emboldened by the success with genus 1 curves, let us move on to the much more subtle case of genus 2.
The generic Sato--Tate group for a genus 2 curve over $\mathbb{Q}$ is $\mathrm{USp}(4)$, which occurs in the case of trivial endomorphism ring. The dominance of the generic case is reflected in the LMFDB, in which 63107 out of 66158 genus 2 curves over $\mathbb{Q}$ have this Sato--Tate group \cite[Genus 2 curves over $\mathbb{Q}$]{lmfdb}. Again, the good Euler factors depend only on the isogeny class. Unlike with elliptic curves, there is no option to ask the LMFDB for one curve per isogeny class. On the other hand, the database has 65534 classes and so over 99\% have a unique representative. With this in mind, we simply accept the redundancy. Using the LMFDB data, we perform the binary classification: Is the Sato--Tate group $\mathrm{USp}(4)$ or not? 

Looking at Eqs.~\eqref{eq.genus2ZL}--\eqref{akp}, we see that the zeta-function for genus 2 curves is governed by an $L$-function numerator which is a degree 4 palindromic polynomial.
Hence, there are two non-trivial (normalized) coefficients, $(a_{1,p}, \ a_{2,p})$ of the Euler factors.
Using SAGE \cite{sage}, we calculate the zeta function of a curve for all first 200 primes $p$ excluding $2$ (i.e., $p < 1230$) which is always bad.

Thus, we can establish the following dataset:
\begin{equation}
\cD := \{ (a_{1,p}, \ a_{2,p})_{2 <p < 1230} \longrightarrow \mbox{ yes/no for } \USp(4) \}
\end{equation}
As mentioned in the opening paragraph, the vast majority are the generic full $\USp(4)$, so we need to down-sample in order to not bias a classifier.
Thus we randomly select 3000 of the $\USp(4)$ cases and combine that with the non-$\USp(4)$ cases (which, from above, is 
$66158- 63107 = 3051$; actually, 2440 of these non-$\USp(4)$ cases belong to the Sato--Tate group $G_{3,3} \cong \SU(2) \times \SU(2)$). On this balanced dataset $\tilde{\cD} \subset \cD$, we again perform cross-validation by taking $20\%$ training, and validating on the remaining $80\%$.
Using a Naive Bayes classifier as the genus 1 case, we here find precision 0.990 and Matthew's phi 0.98, which is excellent.
Again, we have tried other standard classifiers, and we find that nearest neighbours performed similarly, though decision trees were quite a bit worse.

To have an extra confirmation that there is inherent structure in the data.
Let us consider each of the 200 pairs $(a_{1,p}, \ a_{2,p})$ as a point in $\IR^{400}$.
Using principal component analysis (q.v., \cite{ML}), by projecting this point cloud of data from $\IR^{400}$ to $\IR^{2}$, as shown in part (a) of Figure \ref{f:genus2-ap}, we can see that the $\USp(4)$ (marked as 1) and non-$\USp(4)$ (marked as 0) very neatly separate. 

\begin{figure}[!h!t!b]
\centerline{
(a)
\includegraphics[trim=12mm 0mm 0mm 0mm, clip, width=3in]{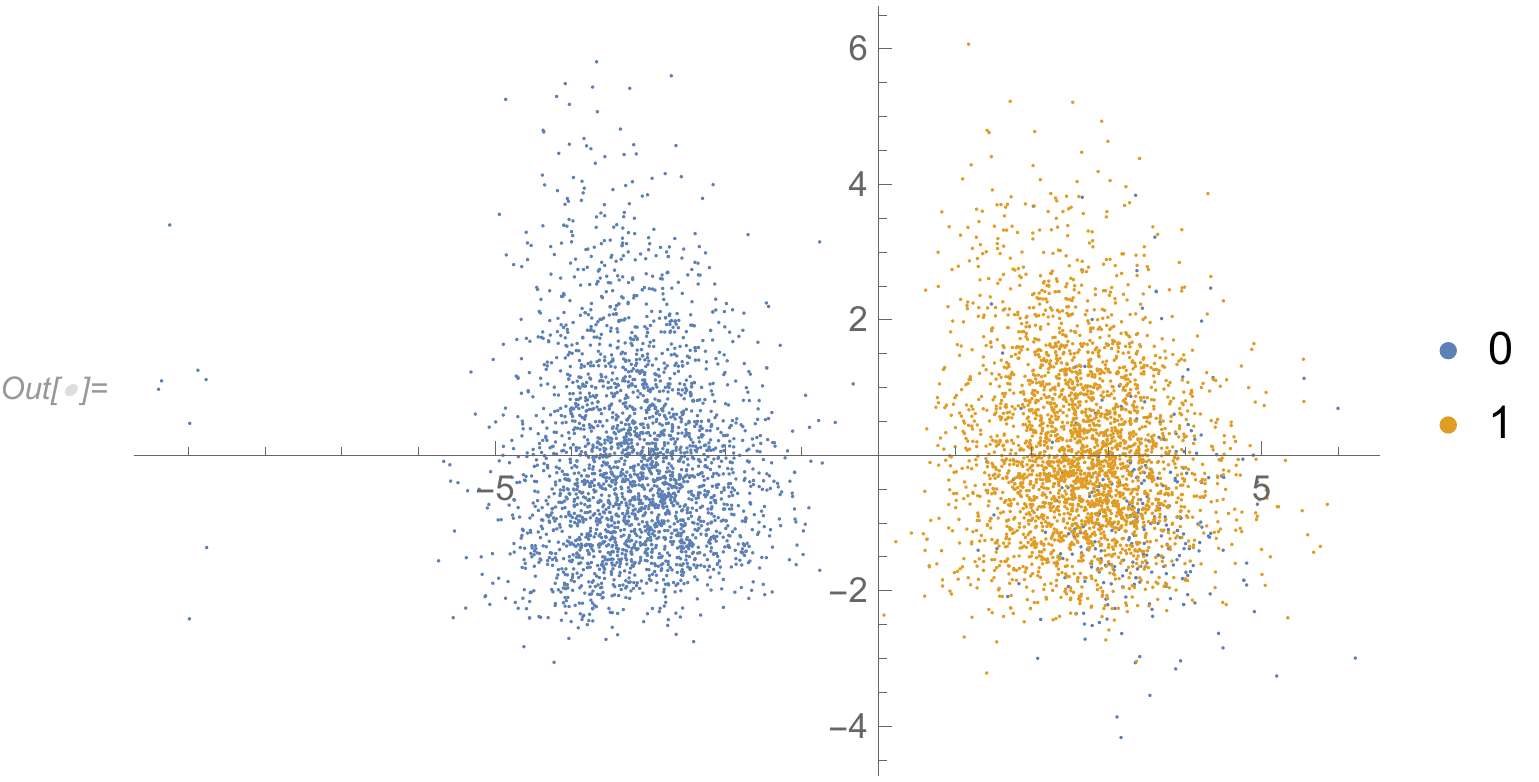}
(b)
\includegraphics[trim=12mm 0mm 0mm 0mm, clip, width=3in]{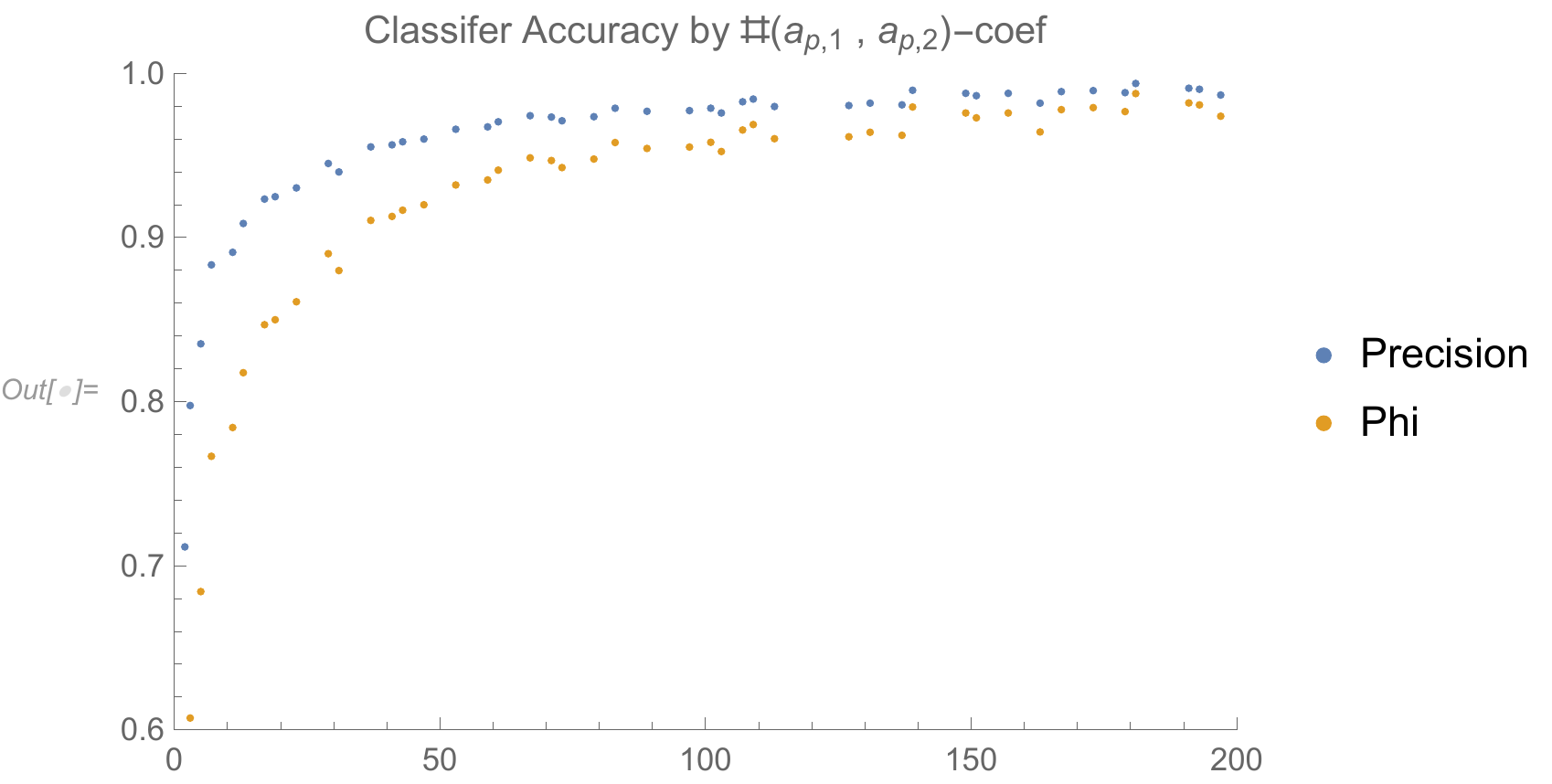}
}
\caption{{\sf {\small
(a)
A principal component analysis (PCA) by projecting the labeled data pairs of coefficients in $\IR^{400}$ corresponding to generic vs. non-generic Sato-Tate to 2-dimension. 
(b)
The precision and confidence (Matthew's phi-coefficient) of the Naive Bayes classifier, for the problem distinguishing the generic Sato-Tate group $\USp(4)$, against the number of $(a_{1,p} , a_{2,p})$ coefficients of the genus 2 hyperelliptic curve supplied to the training.
}}
\label{f:genus2-ap}}
\end{figure}

To get an idea of how effective the training is, we present a gradation of coefficients to the classifier from a single pair (at $p=3$) cumulating to more pairs of $a_p$ coefficients as we go up in primes.
This is drawn in part (b) of Figure \ref{f:genus2-ap}. We see that in the beginning the performance is poor but by the time it has seen around 50 primes, we are already at 0.95 precision. 

\begin{remark}
In a recent paper \cite{Zyw}, D. Zywina shows that one can determine the identity component
of the Sato--Tate group of an abelian variety using just two
$L$-polynomials, though his algorithm does not specify which polynomials 
we need. It would be interesting to consider his result in the perspective of machine-learning.
\end{remark}

\section{Distinguishing non-generic curves using random matrices} \label{sec-finer}

In this section we go beyond the binary classification of the previous section. The Sato--Tate group is a compact Lie group. For genus 2 curves, there are 6 possibilities for its identity component. The non-generic cases occur with decreasing probability, and ultimately the number of occurrences are too small to train the classifier. Worse still, the complete classification of Sato--Tate groups for genus 2 curves over $\mathbb{Q}$ features 34 distinct cases. There is far too little data available on the LMFDB to distinguish these cases by machine learning, for example, only 1 curve on the database has group $D_{6,2}$ \cite[Genus 2 curve 11664.a.11664.1]{lmfdb}.

To circumvent this difficulty, we generate random matrices for training the classifier. The point is that the distribution of the Euler factor coefficients should converge to the distribution of the characteristic polynomial coefficients of random matrices in the Sato--Tate group. This allows us to train classifiers for the non-generic Sato--Tate groups. Still, due to the lack of data, we are unable to verify the classifier's accuracy for curves in certain cases of rare Sato--Tate groups.
Nevertheless, we will see below in several cases where there is sufficient data to verify, the classifier does perform very well.
We keep the notations for the Sato--Tate groups in Definition \ref{def-ST}. 

Specifically, we will do the following, in light of Conjecture \ref{conj-ST}:
\begin{itemize}
\item We fix $k$ different Sato--Tate groups, $ST_{i=1,2,\ldots,k}$, say.
	For each $ST_i$, take 200 random elements within the group (as $4\times4$ matrices) and for each matrix, compute its characteristic polynomial and extract the two non-trivial coefficients $(c_1,c_2)$ as in \eqref{charpoly}.
\item We repeat the above 1000 times. This gives 1000 cases of 200 pairs $(c_1,c_2)$ for each $ST_i$, accordingly labeled.
\item We now train a classifier (Naive Bayes, decision tree, nearest neighbour or otherwise) to this labeled data.
Note that so far, there is no input from number theory or geometry, the classifier has only been fed group-theoretic information: the characteristic polynomial of $ST_i$ matrices.
\item We can now validate the classifier on {\it actual} curve information, viz., for a genus 2 curve from LMFDB, obtain 200 pairs of normalized Euler coefficients $(a_{1,p} , a_{2,p})$ for the first 200 primes $p$. The classifier will then return one of the $k$ labels (categories), which is then compared to the actual Sato-Tate group for the curve.
The precision and confidence for the $k$-category classification is then computed between the predicted and actual.
\end{itemize}
We remark that we are not using moments of the probability distributions. Instead, we are using sample points from the distributions to train a classifier.

\subsection{$N(G_{1,3})$ and $N(G_{3,3})$}

We begin with a binary classification between the non-generic genus 2 Sato--Tate groups $N(G_{1,3})$ and $N(G_{3,3})$. These groups have different identity components. After generating 1000 samples of coefficient pairs for each group, a Bayes classifier can distinguish the corresponding distributions with 100\% accuracy. There are 303 (resp. 144) curves on the LMFDB with group $N(G_{1,3})$ (resp. $N(G_{3,3})$). Given coefficient pairs for the first 200 Euler factors for these curves, the classifier could distinguish the groups with 100\% accuracy. That is, it has completely correctly sorted the 303 vs.~144 genus 2 curves with Sato-Tate group $N(G_{1,3})$ vs.~$N(G_{3,3})$. The running time, again, is less than 1 second on an ordinary laptop, using Mathematica \cite{wolfram}.

To get an idea how many coefficient pairs are needed to efficiently train the classifier, we repeat the above experiment starting with only one pair, and going up gradually. This constitutes a {\it learning curve} where the accuracy and confidence are plotted against the number of pairs seen in the training.
We show this in Part (a) of Figure \ref{f:genus2-ap}.
We see that given only the first coefficient pair, the classifier is useless. At around 10 coefficients its accuracy is already at high 90s, and by 20 or 30 it is all 100\%.

\begin{figure}[!h!t!b]
\centerline{
(a)
\includegraphics[trim=12mm 0mm 0mm 0mm, clip, width=3in]{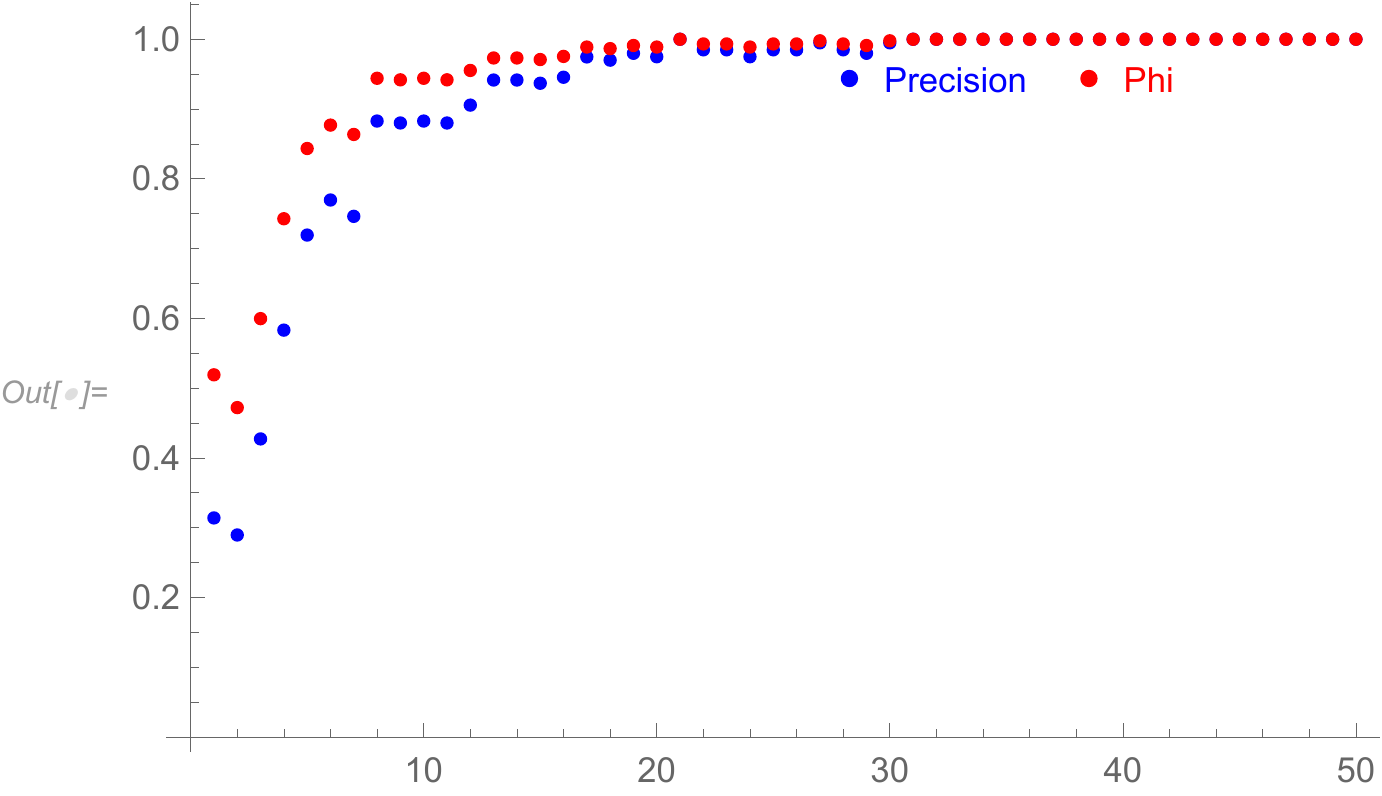}
(b)
\includegraphics[trim=12mm 0mm 0mm 0mm, clip, width=3in]{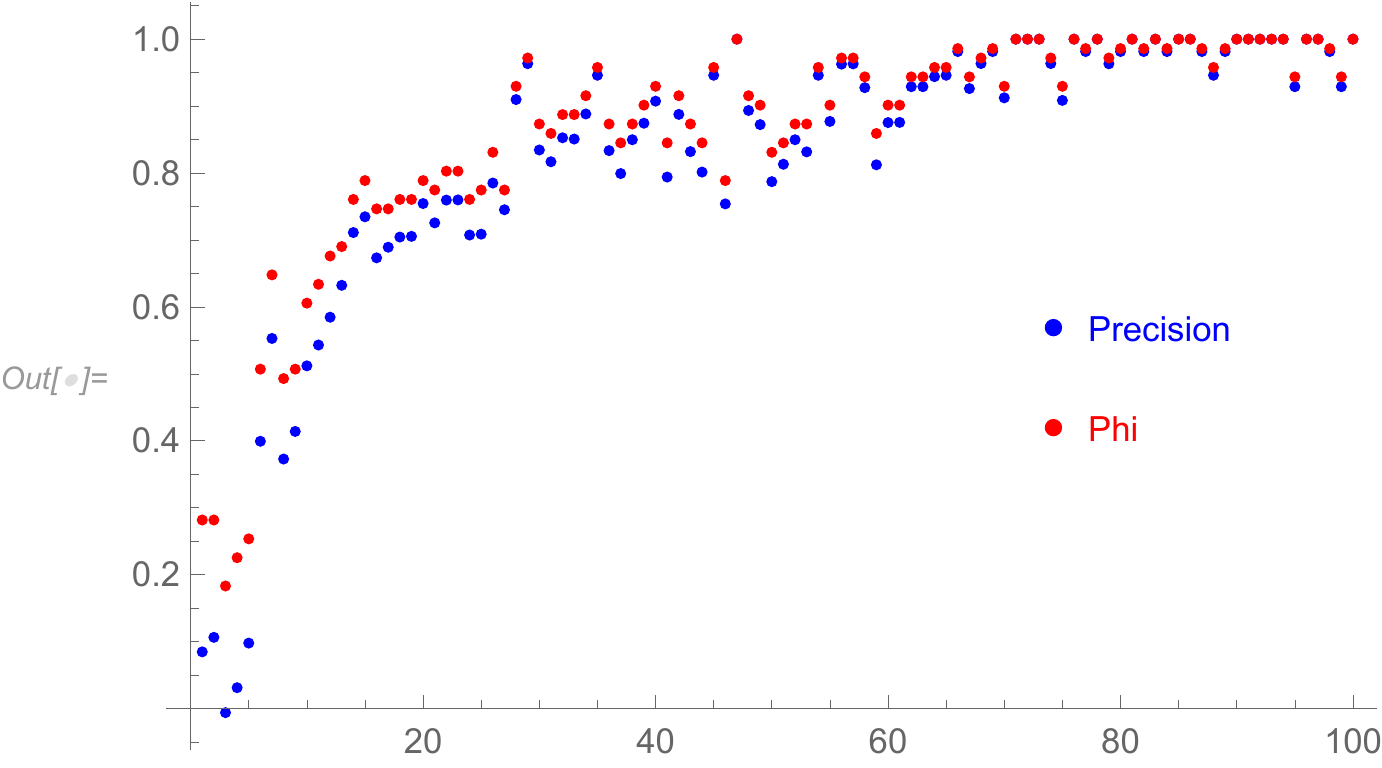}
}
\caption{{\sf {\small
(a)
The precision and confidence  (Matthew's phi-coefficient) of the Naive Bayes classifier, for the problem of distinguishing $N(G_{1,3})$ and $N(G_{3,3})$, against the number of pairs of coefficients $(a_{1,p} , a_{2,p})$ of the genus 2  curve for $p$ up to the value on the X-axis, supplied to the training.
(b)
The same plot, but for the 5-way classifier of the Sato-Tate group $J(E_n)$, $n\in\{1,2,3,4,6\}$.
}}
\label{f:genus2-ap}}
\end{figure}

\subsection{$J(E_n)$, $n\in\{1,2,3,4,6\}$}

We finally attempt a 5-way classification between the non-generic genus 2 Sato--Tate groups $J(E_1)$, $J(E_2)$, $J(E_3)$, $J(E_4)$ and $J(E_6)$. These groups all have the same identity component $\SU(2)$. 
As before, we generate 1000 random samples of 200 coefficient pairs for each of the five $J(E_n)$ groups.
A Naive Bayes classifier is then trained on these.
Upon validating on the actual curve data, of which is there a paucity from LMFDB, a total of 71 cases, we find that the confusion matrix is 
\[
M = {\tiny \left(
\begin{array}{ccccc}
 24. & 0. & 0. & 0. & 0. \\
 0. & 9. & 0. & 0. & 0. \\
 0. & 0. & 3. & 0. & 1. \\
 0. & 0. & 0. & 17. & 0. \\
 0. & 0. & 0. & 0. & 17. \\
\end{array}
\right)},
\] 
which means that only a single case has been mis-classified (the 1 off-diagonal).
The accuracy is  98.59\% and confidence 0.9814.
This is quite impressive for a 5-way classification, in under 1 second.

Again, to get an idea of a learning curve, we show in Part (b) of Figure \ref{f:genus2-ap}, the accuracy and confidence attained by showing an increasing number of coefficients in the training process.
In the beginning the classifier was around 0\% accuracy but by 40-50 coefficient pairs it was getting to almost 100\%. 
All fluctuations are due to the random sampling in the training data.


\section{Conclusion \& Outlook}\label{sec:Outlook}

Let us contrast the efficiency of the Bayes classifier to the established approach for computing the Sato--Tate group. According to \cite[Section~5.2]{FKRS}, with $N=2^{20}$, the first 20 moment statistics for a curve agree with the corresponding moments for the group with an error of $0.1\%$. On the other hand, the best agreement one finds by comparing to other Sato--Tate groups is worse than 40\%. 
The Bayes classifier requires much smaller $N$ (around $2^{10}$), and agrees with the identification via moment sequences. One could use the trained classifier to predict the Sato--Tate groups for arbitrary curves, though we could not verify the accuracy due to a lack of data for rare Sato--Tate groups.

The results in this paper provide convincing evidence that a machine can be trained to learn the Sato--Tate distributions and to classify curves according to their Sato--Tate groups. 
Our approach of using Euler factors is in accordance with the setup of the Langlands program, and we expect that many important objects in number theory can be studied through machine-learning by analyzing data consisting of Euler factors.


{\small
Yang-Hui He {\sf hey@maths.ox.ac.uk} \\
Department of Mathematics, City, University of London, EC1V 0HB, UK;\\
Merton College, University of Oxford, OX14JD, UK;\\
School of Physics, NanKai University, Tianjin, 300071, P.R.~China

Kyu-Hwan Lee {\sf khlee@math.uconn.edu} \\
Department of Mathematics, University of Connecticut, Storrs, CT, 06269-1009, USA

Thomas Oliver {\sf Thomas.Oliver@nottingham.ac.uk} \\
School of Mathematical Sciences, University of Nottingham, University Park, \\
Nottingham, NG7 2QL, UK
}
\end{document}